\documentclass[12pt]{article}

\usepackage{amsmath}
\usepackage{amsthm}
\usepackage{amsfonts}
\usepackage{amssymb}
\usepackage{amscd}

\usepackage{bezier}
\usepackage{enumerate}

\usepackage[dvips]{graphicx}

\theoremstyle{plain}
\newtheorem{maintheorem}{Theorem}

\newtheorem{mainclly}[maintheorem]{Corollary}
\swapnumbers
\newtheorem{thm}{Theorem}[section]
\newtheorem{prop}[thm]{Proposition}

\newtheorem{lemma}[thm]{Lemma}

\newtheorem*{claim}{Claim}

\newcommand{\R}{\mathbb{R}}

\newcommand{\N}{\mathbb{N}}
\newcommand{\nat}{\mathbb{N}}

\newcommand{\cl}{\operatorname{Cl}}

\newcommand{\per}{\operatorname{Per}}

\newcommand{\Diff}{\operatorname{Diff}}

\newcommand{\crit}{\operatorname{Crit}}

\newcommand{\sing}{\operatorname{Sing}}

\title{Homoclinic classes and
finitude of attractors for
vector fields on $n$-manifolds }

\author{C.~M.~Carballo and C.~A.~Morales\thanks{
The first author was supported by CNPq.
The second author was supported by FAPERJ, CNPq,
and PRONEX/Dyn. Sys.}}

\date{March 6, 2001}

\begin{document}

\maketitle

\begin{abstract}
A {\em homoclinic class} of a vector field
is the closure of the transverse homoclinic
orbits associated to a hyperbolic
periodic orbit. 
An {\em attractor} (a {\em repeller})
is a transitive set to which converges
every positive (negative) nearby orbit.
We show that a generic $C^1$ vector field on 
a closed $n$-manifold has
either infinitely many homoclinic classes
or a finite collection of attractors (repellers) 
whose basins form an open-dense set.
This result gives an approach to a conjecture 
by Palis.
We also prove the existence of a locally
residual subset of $C^1$ vector fields
on a $5$-manifold having
finitely many attractors and repellers 
but infinitely many homoclinic classes.
\end{abstract}

\section{Introduction}
\label{s-intro}

In this paper we prove
that a generic $C^1$ vector field
on a closed $n$-manifold has 
either infinitely many homoclinic classes
or a finite collection of attractors (repellers) 
whose basins form an open-dense set.
This result gives an approach to a conjecture 
by Palis in the $C^1$ interior of the set of
vector fields
having finitely many homoclinic classes.
We also prove the existence of a locally
residual subset of $C^1$ vector fields
on a $5$-manifold having
finitely many attractors and repellers but 
infinitely many homoclinic classes.
This proves that the finitude of
attractors and repellers
does not (not even generically) imply
the existence of a (finite)
spectral decomposition
formed by transitive sets.

Before stating our results precisely, we
introduce some definitions and present the
motivations for this work.
The set of $C^1$ vector fields on a closed manifold 
$M$ endowed with the $C^1$ topology
is a Baire topological space that we denote
$\mathcal{X}^1(M)$.
The flow generated by $X\in \mathcal{X}^1(M)$ 
is denoted by $X^t$. 
A subset of a topological space is
{\em residual} if it includes a set that is a countable
intersection of open-dense sets.
Given an open set ${\cal U}$ of $\mathcal{X}^1(M)$
and a property $P$, we say that
{\em a generic vector field in ${\cal U}$ satisfies $P$}
if there is a residual subset ${\cal R}$ of ${\cal U}$ 
such that every element of ${\cal R}$
satisfies $P$.

An invariant set of $X\in \mathcal{X}^1(M)$ is
{\em transitive} if it is the $\omega$-limit
set of one of its orbits.
Recall that the {\em $\omega$-limit set} 
$\omega(x)$ of a point $x$ in $M$ is 
the set of all accumulation
points of the positive orbit of $x$.
Similarly, the {\em $\alpha$-limit set}
$\alpha(x)$ of $x$ is 
the set of all accumulation 
points of the negative orbit of $x$.
The {\em nonwandering} set of $X$,
denoted by $\Omega(X)$, is the set
of $x\in M$ such that
for every neighborhood $U$ of
$x$ and $T>0$ there
is $t>T$ such that
$X^t(U)\cap U\neq\emptyset$. 

An {\em attractor} of $X$ 
is a transitive set $\Lambda$ that 
has a neighborhood $U$ such that
$X^t(U)\subseteq U$, for $t> 0$, and
$\cap_{t\geq 0} X^t(U)= \Lambda$.
A {\em repeller} is an attractor
for the time-reversed vector field
$-X$. 
The {\em basin} of an attractor
$\Lambda$ is the set of
all points $x\in M$ such that
$\omega(x)\subseteq \Lambda$.
Similarly, the basin of a repeller
$\Lambda$ is the set of all
points $x\in M$ such that
$\alpha(x)\subseteq \Lambda$.
If $\Lambda$ is an attractor or a repeller
we denote its basin by $B(\Lambda)$.

A periodic orbit of $X$ is hyperbolic if
its associated Poincar\'e map has no 
eigenvalues with modulus one.
Similarly, a singularity 
$\sigma$ of $X$ is hyperbolic if
the derivative $DX(\sigma)$ has no pure
imaginary eigenvalues.
A {\em sink} is an attracting hyperbolic
closed orbit and a {\em source}
is a sink for $-X$.

Given a hyperbolic periodic orbit
$O$ of a vector field $X$, the sets
\[
W^s(O)= \{x\in M: X^t(x)\to O,\text{ as } t\to \infty\}
\]
and
\[
W^u(O)= \{x\in M: X^t(x)\to O,\text{ as } t\to -\infty\}
\]
are $C^1$ immersed submanifolds of $M$.
A {\em homoclinic class} of $X$ is the closure of the
points of transverse intersection between
$W^s(O)$ and $W^u(O)$, for some hyperbolic
periodic orbit $O$ of $X$.
Homoclinic classes are transitive sets \cite{Ne78}.
Similar definitions and facts hold for diffeomorphisms.

It follows from the definition of homoclinic class
that both attracting and repelling hyperbolic periodic
orbits are homoclinic classes.
So, the finitude of homoclinic classes implies
the finitude of sinks and sources.
Later on we will see (Theorem \ref{thm-main})
that the finitude of homoclinic classes generically
implies the finitude of attractors as well.

A central motivation for this work
is the following theorem by Ma\~n\'e:
a generic $C^1$ surface diffeomorphism 
either has infinitely many sinks or sources 
or it is Axiom A without cycles 
\cite{Man82}.
In particular, if a generic $C^1$ diffeomorphism
has only finitely many sinks and sources, then
it has a finite collections of attractors whose 
basins form an open-dense set and a finite collection 
of repellers whose basins form an open-dense set.
A recent related result \cite{MPa8} implies that
the same conclusion is valid for $C^1$ vector fields
on $3$-manifolds.
Our first result gives a (generically) sufficient
condition for the finitude of attractors for vector 
fields on higher dimensional manifolds.
From now on, let $M$ be a closed $n$-manifold, 
$n\geq 3$.

\begin{maintheorem}
\label{thm-main}
A generic vector field in $\mathcal{X}^1(M)$
has either infinitely many homoclinic classes 
or else a finite collection of attractors
whose basins form an open-dense set 
and a finite collection of repellers 
whose basins form an open-dense set.
\end{maintheorem}

As a consequence of this, we have the following
corollary dealing with the global conjecture
on the finitude of attractors and their metric 
stability by Palis \cite{Pa00}.
Let $\mathcal{HC}^1(M)$ be 
the $C^1$ interior of the set of 
$X\in \mathcal{X}^1(M)$ such that 
$X$ has finitely many homoclinic classes.

\begin{mainclly}
\label{cll-pc}
If $X\in \mathcal{HC}^1(M)$
and $\Lambda$ is an attractor,
then there is a negihborhood 
${\cal U}\subseteq \mathcal{HC}^1(M)$
of $X$ such that a generic vector field in 
${\cal U}$ (i.e. a generic perturbation of $X$)
has a finite collection of attractors 
whose basins contain an open-dense subset 
of $B(\Lambda)$.
A generic vector field in $\mathcal{HC}^1(M)$
has a finite collection of attractors whose basins
form an open-dense set.
\end{mainclly}

The result in \cite{MPa8} mentioned above shows that 
if a generic $C^1$ vector field on a closed $3$-manifold
has finitely many attractors and repellers, then
it also has finitely many homoclinic classes.
This is not true in higher dimensions by the following
example.

\begin{maintheorem}
\label{thA}
There are a closed $5$-manifold $M$ and 
an open subset ${\cal U}$ of $\mathcal{X}^1(M)$
such that a generic vector field in ${\cal U}$ has
two attractors
(whose basins form an open-dense set),
a unique repeller (with open-dense basin),
and infinitely many homoclinic classes.
\end{maintheorem}

We remark that both 
Theorem \ref{thm-main} and Corollary \ref{cll-pc} 
hold for diffeomorphisms.
Theorem \ref{thA} also holds for diffeomorphisms
on $n$-manifolds, $n\geq 4$, as we will see later
(Proposition \ref{ex1}).

\section{Proof of 
Theorem \ref{thm-main} and Corollary \ref{cll-pc}}
\label{s-pf}

We start by summarizing the generic properties of 
$C^1$ vector fields that we will use and we establish a 
preliminary lemma that gives a sufficient condition for
an invariant set to be an attracting set.

The union of the closed orbits of a vector field
$X$ is denoted by $\crit(X)$.
Pugh's General Density Theorem \cite{Pu67}
defines a residual subset $\mathcal{P}$ of
$C^1$ vector fields which have the following
properties.
If $X\in \mathcal{P}$, then
all its closed orbits are hyperbolic,
the stable and unstable manifolds of
the closed orbits intersect transversely, and
$\Omega(X)= \cl(\crit(X))$.
In particular, $X$ has a finite number 
of singularities.

We denote by $\mathcal{Q}$ the residual subset
of $C^1$ vector fields introduced in \cite{CMP1}.
If $X\in \mathcal{Q}$, then its homoclinic classes
either coincide or are disjoint.

A compact set $\Lambda$ of a $C^1$
vector field $X$ is {\em isolated} if it has a
neighborhood $U$ (called {\em isolating block}) 
such that
$\cap_{t\in \R}X^t(U)= \Lambda$.
An {\em attracting set} is an isolated
set having an isolating block $U$ 
such that $X^t(U)\subset U$ for
every $t>0$.
A compact set $\Lambda$ is
{\em Lyapunov stable} for $X$ if 
for every neighborhood $U$ of $\Lambda$ 
there is another neighborhood
$V\subseteq U$ such that $X^t(V)
\subset U$ for every $t>0$.

Given a vector field
$X$ on a closed $n$-manifold
$M$ we denote
by $R_X$ the set of all $x\in M$
such that $\omega(x)$ is Lyapunov stable 
for $X$.
We will use the following result.

\begin{thm}(\cite{MPa6})
\label{mpa}
The set
\[
\mathcal{S}= \{X\in \mathcal{X}^1(M):
R_X \mbox{ is residual in }M\}
\]
is
residual in $\mathcal{X}^1(M)$.
\end{thm}

The attractors of a vector field are 
transitive isolated Lyapunov stable sets.
The lemma below implies the converse.
It seems to be well known and we prove 
it here for completeness.

\begin{lemma}
\label{lem-csas}
Isolated Lyapunov stable sets are
attracting sets.
\end{lemma}

\begin{proof}
Let $\Lambda$ be an isolated Lyapunov stable set
of a vector field $X$.
We first observe that there is an
open set sequence
$U_n$ such that $ U_{n+1}\subseteq U_n$,
$X^t(U_n)\subset U_n$,
for all $n$, and
$\Lambda= \cap_{n\in \N} U_n$.
To see this, let
$B_n$ be the ball of radius $1/n$
around $\Lambda$, i.e. 
$B_n= \{ x\in M: d(x,\Lambda) < 1/n \}$.
By the Lyapunov stability of $\Lambda$, for each $n$, 
there is an open neighborhood $V_n\subseteq B_n$ of 
$\Lambda$ such that
$X^t(V_n)\subseteq B_n$ for $t\geq 0$.
Setting $U_n= \cap_{t\geq 0} X^t(V_n)$ we are done.

Now, let $U$  be an isolating block for $\Lambda$ and 
$n$ be such that $U_n\subseteq U$.
Then, $X^t(U_n)\subset X^t(U)$
for every $t\in \R$. 
So,
\begin{equation}
\label{e1}
\cap_{t\in \R}X^t(U_n)\subseteq \Lambda
\end{equation}
since $U$ is an isolating block of $\Lambda$.

On the other hand,
as $U_n$ is forward invariant, we have that
$U_n\subseteq X^t(U_n)$ for  $t\leq 0$.
Then,
$ U_n\subset \cap_{t\leq 0}X^t(U_n) $ and so
\[
\cap_{t\geq 0}X^t(U_n)
=U_n\cap\left(\cap_{t>0}X^t(U_n)\right)
\subset\left(\cap_{t\leq 0}X^t(U_n)\right)
\cap
\left(\cap_{t>0}X^t(U_n)\right)
=\cap_{t\in \R}X^t(U_n).
\]
Applying (\ref{e1}) we obtain
$\cap_{t\geq 0}X^t(U_n)\subset \Lambda$.
But $\Lambda\subset U_n$ and $\Lambda$ is invariant.
Then, $\Lambda\subset X^t(U_n)$ for every $t\geq 0$
and so $\Lambda\subset \cap_{t\geq 0}X^t(U_n)$.
This proves
$\cap_{t\geq 0}X^t(U_n)= \Lambda$.
Recalling that $U_n$ is a neighborhood
of $\Lambda$, we see that $\Lambda$ is an
attracting set since it is compact.
This proves the lemma.
\end{proof}

{\flushleft{\bf Proof of Theorem~\ref{thm-main}: }}
We define a residual subset $\mathcal{R}$ of 
$\mathcal{X}^1(M)$ by
$\mathcal{R}= 
\mathcal{P}\cap \mathcal{Q}\cap \mathcal{S}$ 
and consider $X\in \mathcal{R}$ so that 
our hypothesis implies that $X$ satisfies all the
generic properties above.

If $X$ has infinitely many homoclinic classes, then we
are done.
If not, let $H_1,\dots, H_m$ be the collection of
all homoclinic classes of $X$.
By hypothesis, Pugh's General Density Theorem applies
to $X$; hence $X$ has finitely many singularities and 
$\Omega(X)= \cl(\crit(X))$.
Let $\sing(X)$ be the set of singularities of $X$ and
$\per(X)$ be the set of periodic points of $X$.
Then, we have that $\Omega(X)= \sing(X)\cup \cl(\per(X))$.
On the other hand, by the definition of homoclinic classes,
it is clear that $\cl(\per(X))\subseteq \cup_{i=1}^m H_i$
(the closure of the periodic points is contained in the
closure of the union of the collection of homoclinic 
classes).
From this, 
$\Omega(X)\subseteq \sing(X)\cup (\cup_{i=1}^m H_i)$.
Observe that we also have the other inclusion because 
each homoclinic class is included in the nonwandering 
set of $X$.

Then,
\[
\Omega(X)= \sing(X)\cup (\cup_{i=1}^m H_i).
\]
Let $\sigma_1,\dots, \sigma_l$ be the 
singularities that are not contained in any homoclinic
class; we have that
\[
\Omega(X)= 
\{\sigma_1\}\cup\dots\cup\{\sigma_l\}\cup
H_1\cup\dots\cup H_m.
\]

Observe that, by hypothesis, $X$ satisfies the generic
properties in \cite{CMP1}. In particular,
the collection $\{H_i:i=1,\cdots , m\}$ is disjoint.
The last equality and the fact that
$\{H_i:i=1,\cdots , m\}$ is disjoint
imply that
every $H_i$ is {\em $\Omega$-isolated},
i.e. $\Omega(X)\setminus H_i$ is a closed set.
By \cite{CMP1} we conclude that
every $H_i$ is an isolated set.

Let $\Lambda$ be an attractor of $X$ and $x\in M$ such that
its $\omega$-limit set $\omega(x)= \Lambda$.
In particular, $\omega(x)$ is Lyapunov stable.
We claim that $\Lambda$ is one of the sets in the decomposition. 
Indeed, we have the following:

\begin{claim}
Let $X$ as above and $x\in M$.
If $\omega(x)$ is Lyapunov stable, then 
$\omega(x)$ is either one of the $\sigma_i$'s 
or one of the $H_i$'s.
\end{claim}

\begin{proof}
Clearly $\omega(x)\subset \Omega(X)$.
By assumption we have
$\Omega(X)= 
\{\sigma_1\}\cup\dots\cup \{\sigma_l\}\cup
H_1\cup\dots\cup H_m$ 
(disjoint union).
Suppose that
$\omega(x)\cap H_i\neq \emptyset$ 
for some $i$.
Then, $\omega(x)\subseteq H_i$ since
$\omega(x)\subseteq \Omega(x)$ is connected and 
the above union is disjoint.
As $\omega(x)$ is Lyapunov stable
and $H_i$ is transitive,
we conclude that
$H_i\subseteq \omega(x)$.
We conclude that $\omega(x)= H_i$.

Next suppose that
$\omega(x)\cap \{\sigma_i \}\neq \emptyset$ 
for some $i$. 
Then, $\omega(x)\subseteq \{\sigma_i\}$
by connectedness once more.
We conclude that $\omega(x)= \{\sigma_i \}$.
To finish we note that
$\omega(x)$ must intersect either 
$\{\sigma_i\}$ or $H_i$ (for some $i$) since
$\omega(x)\subseteq \Omega(X)$
and
$\Omega(X)= 
\{\sigma_1\}\cup\dots\cup \{\sigma_l\}\cup
H_1\cup\dots\cup H_m$.
This proves the claim.
\end{proof}

The Claim implies that $X$ has only finitely
many attractors 
(because any of them is either $\sigma_i$ or
$H_i$ for some $i$).
Let us prove that the union of the basins of the
attractors of $X$ is an open-dense set of $M$.
By hypothesis $X\in \mathcal{S}$.
So, by Theorem \ref{mpa}, the set
$R_X$ of
$x\in M$ such that $\omega(x)$ is Lyapunov stable is
residual in $M$.
Denote
$B(X)= \cup\{B(\Lambda): \Lambda \text{ is an attractor of } X \}$.
We will prove that $R_X\subseteq B(X)$. 
Indeed, choose $x\in R_X$. 
Then,
$\omega(x)$ is
Lyapunov stable.
The Claim implies that $\omega(x)$ is
either $\sigma_i$ or $H_i$ for some $i$.
If $\omega(x)= \{\sigma_i\}$ for some $i$,
then we have that
$\sigma_i$ is a sink since
$\sigma_i$ is hyperbolic and
$\omega(x)$ is Lyapunov stable.
In particular, $\omega(x)$ is an attractor
of $X$.
If $\omega(x)=H_i$ for some $i$ we have that
$H_i$ is Lyapunov stable.
As $H_i$ is isolated we can apply
Lemma~\ref{lem-csas}
to obtain that $H_i$ is an attracting set.
As $H_i$ is transitive we conclude that
$H_i$ (and hence $\omega(x)$) is
an attractor.
So,
$\omega(x)$ is an attractor
in this case as well.
This proves that $R_X\subseteq B(X)$.
As $R_X$ is residual we conclude that
$B(X)$ is dense.
As $B(X)$ is obviously open we conclude
that $B(X)$ is an open-dense set of $M$.

We have proved that a generic 
$X\in \mathcal{X}^1(M)$ has 
either infinitely many homoclinic classes 
or a finite collection of
attractors with open-dense basin.
The same argument applied to the 
time-reversed vector field $-X$ shows that
a generic $X\in \mathcal{X}^1(M)$ satisfies
a similar conclusion for repellers 
(instead of attractors).
The proof follows since the intersection of
residual sets is residual as well.
\qed
\bigskip

{\flushleft{\bf Proof of Corollary~\ref{cll-pc}: }}
The second part of Corollary \ref{cll-pc} follows 
directly from Theorem \ref{thm-main}.
To prove the first part, 
we let $X\in \mathcal{HC}^1(M)$ and $\Lambda$ be
an attractor of $X$.
Let $Y$ be a generic $C^1$ perturbation of $X$ so that
we can assume that $Y\in \mathcal{R}\cap \mathcal{HC}^1(M)$.
By Theorem~\ref{thm-main} we have that $Y$ has
finitely many attractors whose basins of
attraction form an open-dense set of $M$.
As before, we let $B(Y)$ denote the union of the basins
of the attractors of $Y$.
Now, let $\Lambda_1,\dots,\Lambda_k$ be the attractors of $Y$
whose basins meet the basin of $\Lambda$; i.e. such that
$B(\Lambda_i)\cap B(\Lambda)\neq \emptyset$, for $i=1,\dots,k$.
We must show that $(\cup_{i=1}^k B(\Lambda_i))\cap B(\Lambda)$ 
is an open-dense subset of $B(\Lambda)$.
It is clear that it is open.
It is also clear that it is dense because
$(\cup_{i=1}^k B(\Lambda_i))\cap B(\Lambda)\supseteq 
B(Y)\cap B(\Lambda)$ 
and we know that $B(Y)$ is dense in $M$.
This proves the result.
\qed

\section{Proof of Theorem \ref{thA}}

In this section we prove Theorem \ref{thA}.
The proof is a simple application of the
normally hyperbolic theory and \cite{BD99}.
The vector field (and the manifold) in the proposition will
be obtained by suspending a suitable
$C^\infty$ diffeomorphism $F$ in 
the $4$-sphere $S^4$
exhibiting a normally expanding
submanifold diffeomorphic to 
the $3$-sphere $S^3$.
In such submanifold $F$ has
two saddle-type fixed points
(with complex eigenvalues)
persistently connected in the sense of 
\cite{BD99}.
Appart from this, $F$
has two sinks attracting every positive
orbit outside $S^3$.
It would be interesting to prove
the existence of a closed manifold of 
dimension $n\leq 4$ satisfying
the conclusion of Theorem \ref{thA}
(see also \cite{Ab01}).

Throughout, $M$ is a closed manifold
and $\Diff^r(M)$
is the space of $C^r$ diffeomorphisms of $M$,
$r\geq 1$.
The proof of Theorem \ref{thA}
follows by suspending the following 
example.

\begin{prop}
\label{ex1}
There are $F\in \Diff^1(S^4)$,
a neighborhood 
${\cal U}\subseteq \Diff^1(S^4)$ of $F$ and
a residual subset
${\cal R}$ of ${\cal U}$ such that:
\begin{enumerate}
\item
if $G\in {\cal U}$, then
$G$ has two attractors (whose basin form
an open-dense set) and a repeller
(with open-dense basin), and
\item
if $G\in {\cal R}$, then
$G$ has infinitely many homoclinic classes.
\end{enumerate}
\end{prop}

To prove this proposition we need to 
recall some concepts of normally hyperbolic theory.
Given $F\in \Diff^1(M)$ and a closed
invariant submanifold
$N\subset M$ of $F$,
we say that
$N$ is {\em normally expanding}
for $F$ if there are constants $\sigma>1$
and a continuous invariant splitting
$TM=E^u_N\oplus TN$
over $N$ such that
\[
\| DF^{-n}(x)/E^u_x \|
\leq \sigma^{-n},
\]
and
\[
\| DF^{-n}(x)/E^u_x \|
\cdot\| DF^n(F^{-n}(x))/TN_{F^{-n}(x)} \|
\leq \sigma^{-n}
\]
for every $n\in \nat$ and $x\in N$.

\begin{lemma}
\label{exemplo-1}
Let $M$ be a closed manifold,
$F\in \Diff^1(M)$, and $N$ be a normally
expanding manifold of $F$.
Suppose that
there is a neighborhood 
${\cal U}\subseteq \Diff^1(N)$ of
$F/N$ and a residual subset ${\cal R}$ of
${\cal U}$ such that:
\begin{enumerate}
\item
if $g\in {\cal U}$, then $g$ has
a unique repeller, and
\item
if $g\in {\cal R}$, then
$g$ has infinitely many
sinks.
\end{enumerate}
Then, there is a neighborhood ${\cal U}'
\subseteq \Diff^1(M)$ of $F$
and a residual subset ${\cal R}'$ of
${\cal U}'$ such that:
\begin{itemize}
\item
if $G\in {\cal U}'$, then 
$G$ has a unique repeller close to $N$, and
\item
if $G\in {\cal R}'$, then
$G$ has infinitely many homoclinic classes.
\end{itemize}
\end{lemma}

\begin{proof}
Note that $N$ is a normally hyperbolic manifold
of $F$ (see \cite{PT77}).
It follows that there is a neighborhood
${\cal U}'_0$
of $F$ where the normally expanding manifold
$N$ has a continuation $N_G$,
for every $G\in {\cal U}'_0$.
Moreover, for every $G$ which is $C^1$ close to $F$
there is a $C^1$ diffeomorphism $h_G: N_G\to N$.
$h_G$ varies continuously with $G$ and
$h_F: N\to N$ is the identity map.

Using $h_G$ we define
$g_G\in \Diff^1(N)$
by
\[
g_G= h_G\circ (G/N_G)\circ h_G^{-1}.
\]
The map
$\Phi: G\mapsto g_G$ is continuous in ${\cal U}'_0$
(see the remark in Appendix 1 of \cite{PT77})
and $\Phi(F)= F/N$.

As ${\cal R}$ is residual in ${\cal U}$ we
have that
${\cal R}=\cup_n{\cal U}_n$
for some open-dense set
sequence ${\cal U}_n$ of ${\cal U}$.
Define
${\cal U}'=\Phi^{-1}({\cal U})$
and ${\cal U}_n'=\Phi^{-1}({\cal U}_n)$
for all $n$.
Clearly both ${\cal U}'$ and ${\cal U}'_n$ are
open sets, and ${\cal U}_n'\subseteq {\cal U}'$.
It is also clear that every $G\in {\cal U}'$
has a unique repeller close to $N$.

We claim that ${\cal U}'_n$ is dense
in ${\cal U}'$ for every $n$.
Indeed, fix $n\in \nat$ and
$G\in {\cal U}'$.
Thus, $\Phi(G)= g_G\in {\cal U}$.
The denseness of ${\cal U}_n$ in ${\cal U}$
allows us to choose a sequence
$g_k\in {\cal U}_n$ such that 
$g_k\to g_G$, as $k\to \infty$.
Define $\widehat{g}_k: N_G\to N_G$
by
\[
\widehat{g}_k=
h_G^{-1}\circ g_k\circ h_G.
\]
It follows that $\widehat{g}_k\to G/N_G$, 
as $k\to \infty$.
Then, for every $k$ large
there is an extension
$G_k\in \Diff^1(M)$ of
$\widehat{g}_k$ to $M$ so that the sequence
$G_k$ converges to $G$
(see, for instance, \cite[p. 8]{PM82} or 
the proof of Theorem 3 in \cite[p. 107]{Lim99}).
Indeed, there is a sequence
$G_k\in \Diff^1(M)$ such that $G_k\to G$, 
as $k\to \infty$, 
$G_k(N_G)= N_G$, and 
$G_k/N_G= \widehat{g}_k$
for all $k$ large.
In particular, $G_k\in {\cal U}'$.
The uniqueness of the continuation of $N$ 
for $G$ close to $F$ then implies
that
$h_{G_k}= h_G$,
$N_{G_k}= N_G$, and
$\Phi(G_k)=
h_G\circ (G_k/N_G)\circ h_G^{-1}=
g_k\in {\cal U}_n'$.
This proves that
$G_k\in {\cal U}_n'$.
The claim follows since
$G_k\to G$, as $k\to \infty$.

By the claim we have that
${\cal R}'= \cap_n{\cal U}_n'$ is
residual in ${\cal U}'$.
Observe that if $G\in {\cal R}'$,
then $\Phi(G)\in {\cal U}_n$,
for all $n$, yielding
$\Phi(G)\in {\cal R}$.
We conclude that $\Phi(G)$ has
infinitely many sinks.
As $G/N_G$ and $\Phi(G)$ are conjugated
we conclude that $G/N_G$ has infinitely many sinks
too.
But $N_G$ is normally expanding, so,
the infinitely many sinks
of $G/N_G$ correspond
to infinitely many saddle-type
hyperbolic periodic
orbits of $G$ whose homoclinic
classes are trivial. 
We conclude that
every $G\in {\cal R}'$ has 
infinitely many homoclinic classes 
and the lemma follows.
\end{proof}

Let $N$ be a closed manifold and
$f\in \Diff^1(N)$.
Following \cite{BD99} we say that
two hyperbolic periodic points $P', Q'$ of
$f$ are
{\em persistently connected}
if there is a neighborhood 
${\cal U}\subseteq \Diff^1(N)$ of
$f$ and a dense subset
${\cal D}$ of ${\cal U}$ such that
for every $g\in {\cal D}$ there is
a transitive set of $g$ containing
the continuations of both $P'$ and $Q'$.

\begin{thm}(\cite[Theorem C]{BD99})
\label{BD}
Let $N$ be a closed $3$-manifold,
$f\in \Diff^1(N)$ and $P', Q'$ be
different hyperbolic fixed points of $f$
which are persistently connected.
Suppose that
$P'$ has a nonreal complex eigenvalue 
with modulus $<1$, that
$Q'$ has a nonreal complex eigenvalue
with modulus $>1$, and that 
the Jacobian of $f$ at
$P'$ is $<1$.
Then there is a neighborhood ${\cal U}$ of
$f$ and a residual subset
${\cal R}$ of ${\cal U}$ such that
if $G\in {\cal R}$,
then $G$ has infinitely many sinks.
\end{thm}

This result is used to prove the following lemma.

\begin{lemma}
\label{exemplo-2}
There are $f_0\in \Diff^1(S^3)$,
a neighborhhod ${\cal U}\subseteq \Diff^1(S^3)$ of $f_0$,
and a residual subset ${\cal R}$ of ${\cal U}$ such that:
\begin{enumerate}
\item
if $g\in {\cal U}$, then $g$ has a unique repeller, and
\item
if $g\in {\cal R}$, then $g$ has infinitely many sinks.
\end{enumerate}
\end{lemma}

\begin{figure}[htv]
\begin{center}
\includegraphics[scale=0.5]{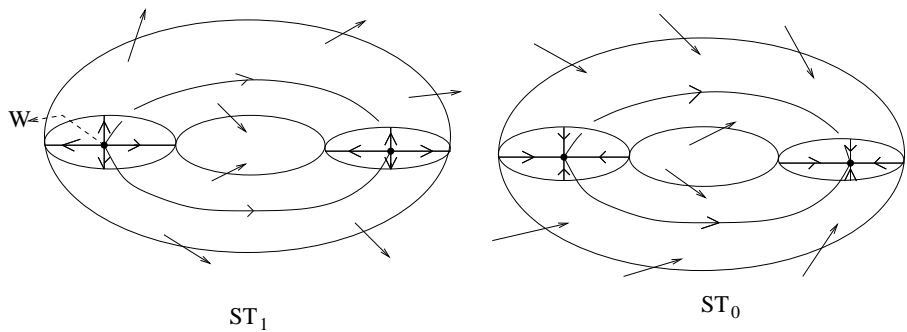}
\end{center}
\caption{The map $g_0$.}
\label{f.1} 
\end{figure}

\begin{figure}[htv]
\begin{center}
\includegraphics[scale=0.5]{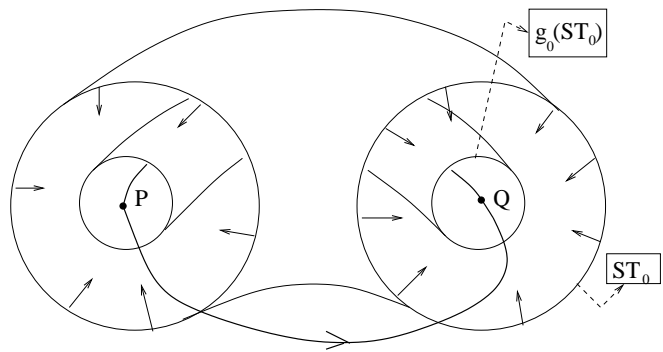}
\end{center}
\caption{Configuration in $ST_0$.}
\label{f.2} 
\end{figure}

\begin{proof}
Consider $S^3$ as
two solid tori $ST_1,ST_0$ glued along
their boundary tori.
Let $g_0\in \Diff^\infty(S^3)$ be
as in Figure \ref{f.1}.
Note that $g_0$ has four fixed points:
two saddles, one sink, and one source, 
denoted by $W$ in Figure \ref{f.1}.
The orbits in $ST_1$ exit $ST_1$ and go to $ST_0$
except for the ones in the middle circle.
We obtain the required diffeomorphism $f_0$
by modifying $g_0$ inside $ST_0$.

First, we observe that
the configuration inside the solid torus
$ST_0$ is as in Figure \ref{f.2}.
In particular, $ST_0$ is contracted
inside it in the way described in Figure \ref{f.2}.
Two fixed points $P, Q$ are indicated
in that figure.

We deform
$g_0(ST_0)$, the image of $ST_0$ by
$g_0$, to obtain $g_1\in \Diff^1(S^3)$
as in Figure \ref{f.3}.
The modification produces
a horseshoe in the meridian disks of $ST_0$
as indicated; observe that
$Q$ and $P$ are fixed by $g_1$.
The modification also
produces two fixed points
$Q', P'$ of $g_1$ described in Figure \ref{f.3}.
Note that $Q'$ (resp. $P'$)
and $P$ (resp. $Q$) are homoclinically related.
We also require the resulting map
$g_1$ to be dissipative
(Jacobian $<1$) in $ST_0$.

\begin{figure}[htv]
\begin{center}
\includegraphics[scale=0.6]{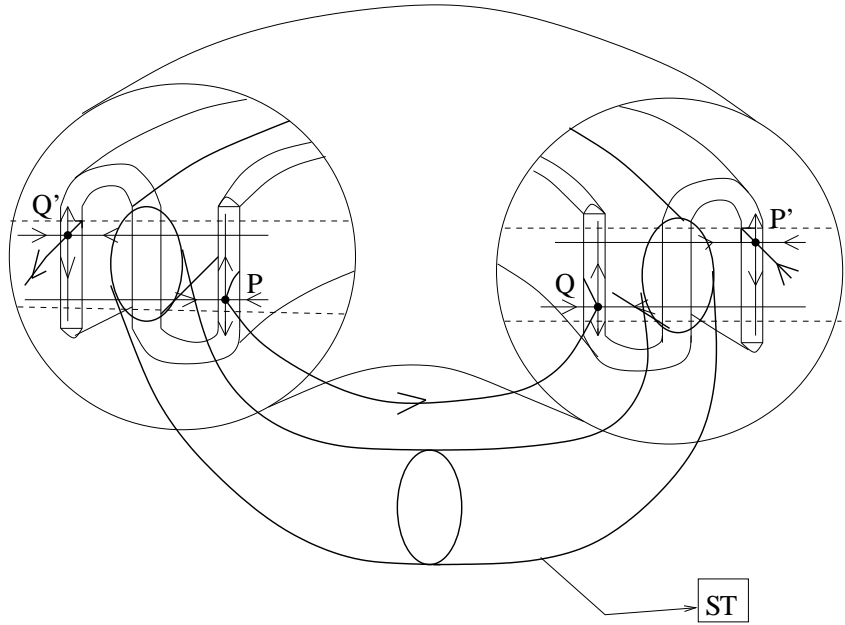}
\end{center}
\caption{Creating a heterodimensional cycle.}
\label{f.3} 
\end{figure}

Next, we modify $g_1$ to obtain
$g_2\in \Diff^1(S^3)$
by performing a rotation supported on
the middle solid torus $ST$
indicated in Figure \ref{f.3}.
This rotation is done to create
a heterodimensional cycle
involving $P,Q$ 
(see \cite{BD99}).
We still have that the map $g_2$ is
dissipative in $ST_0$.

Finally, we modify $g_2$ in order to obtain
$g_3\in \Diff^1(S^3)$ such that
the expanding eigenvalues of
$Q'$ (resp. contracting eigenvalues of $P'$)
become complex conjugated. This can be done
without destroying the heterodimensional cycle
involving $P,Q$.

Unfolding the cycle we create a blender
(``m\'elangeur'') containing $P,Q$ 
(see \cite{BD96}).
The existence of the blender
guarantees that $P'$ and $Q'$ are
persistently connected.
This procedure is similar to the one described
in \cite[p. 149]{BD99}.

Define $f_0=g_3$.
By Theorem \ref{BD}
it follows that there is a neighborhood
${\cal U}$ of $f_0$ and
${\cal R}$ residual in ${\cal U}$ such that
every $g\in {\cal R}$ has infinitely many sinks.
We observe that every $g\in {\cal U}$
has a unique repeller:
the continuation
of the source $W\in ST_1$ of $g_0$
indicated in Figure \ref{f.1}
(note that $W$ is also a source of $f_0$).
This is because any other repeller
of $g$ close to $f_0$ must belong
to the dissipative region $ST_0$
but this is impossible 
by the dissipative condition.
This proves the result.
\end{proof}

{\flushleft{\bf Proof of Proposition \ref{ex1}: }}
Let $f_0\in \Diff^1(S^3)$ be
as in the previous lemma
(Lemma \ref{exemplo-2}).
Let $B$ be a closed $4$-ball in
$\R^4$ centered at $b$.
Note that $\partial B=S^3$.
We extend $f_0\in \Diff^1(S^3)$
to $F_0\in \Diff^1(B)$ such that
$S^3$ is a normally expanding
submanifold of $F_0$
and $b$ is a sink of $F_0$ attracting
every forward orbit of $F_0$ outside $S^3$.
Consider the double manifold $2B=S^4$ and 
the double map $F=2F_0$.
It follows that
$N=S^3$ is a normally expanding
submanifold of $F$
and that $F$ has two sinks
attracting all the forward orbits
outside $N$.
The properties (1) and (2) of
Lemma \ref{exemplo-2}
imply that $F/N=f_0$
satisfies the properties
(1) and (2) of
Lemma \ref{exemplo-1}
for some neighborhood
${\cal U}\subseteq \Diff^1(N)$
and some residual subset ${\cal R}$
of ${\cal U}$.
Then, by Lemma \ref{exemplo-1},
there is a neighborhood
${\cal U}'\subseteq \Diff^1(S^4)$
of $F$ and a residual subset
${\cal R}'$ of ${\cal U}'$ such that
every $G\in {\cal U}'$ has only
one repeller close to
$N$ and every $G\in{\cal R}'$ has
infinitely many homoclinic classes.
We conclude that every
$G\in {\cal U}'$ has a unique
repeller since
any other repeller must be far
from $N$, and so, it must be
the continuation of one of the sinks
of $F$, a contradiction.
Finally we observe that
every $G$ close to $F$ has
two attractors:
the continuation of
the sinks of $F$
(whose basin form an open-dense set of $S^4$),
and a unique
repeller which is
the continuation of the source
$W$ in  Figure \ref{f.1}
(with open-dense basin).
This completes the proof.
\qed
\bigskip

Using the methods of this section
one can prove the
following result.
There is an open set
${\cal U}\subseteq \Diff^1(S^3\times S^1)$
and a residual subset 
${\cal R}\subseteq {\cal U}$
such that:
(1) if $g\in {\cal U}$, 
then $g$ has an
attractor and a repeller
(both with open-dense basin),
and
(2) if $g\in {\cal R}$, 
then $g$ has infinitely
many homoclinic classes.
This allows us to reduce the number of
attractors (repellers) to one,
the minimal possible.

It is interesting to observe that the
diffeomorphisms $F$ in the example of 
Proposition \ref{ex1} have a continuous
dominated splitting 
$TM_{\Omega(F)}= E^s\oplus E^c\oplus E^u$
such that $E^s$ is contracting, $E^u$
is expanding, and either $E^s$ or $E^u$ is
nontrivial.
This shows that the existence of such a
splitting does not imply the existence of
a (finite) spectral decomposition of the 
nonwandering set.

\bigskip

\bibliography{bcc}
\bigskip

\begin{tabular}{ll}
C.~M.~Carballo& C. A. Morales\\
PUC-Rio&
Instituto de Matem\'atica \\
Dto. de Matem\'atica&
Universidade Federal do Rio de Janeiro \\
Rua Marqu\^es de S\~ao Vicente, 225&
C. P. 68.530\\
22453-900, Rio de Janeiro, RJ&
21945-970, Rio de Janeiro, RJ\\
Brazil&
Brazil\\
\verb"carballo@mat.puc-rio.br"&
\verb"morales@impa.br"
\end{tabular}

\end{document}